\documentclass[11pt]{article}
\usepackage{amsmath}
\usepackage{amsfonts}
\usepackage{bm}
\usepackage{cases}
\usepackage{color}
\usepackage{graphicx}

\newtheorem{de}{Definition} 
\newtheorem{thm}{Theorem}
\newtheorem{lem}{Lemma}
\newtheorem{rem}{Remark}

\newcommand{\R}{\mathbb{R}}
\newcommand{\C}{\mathbb{C}}

\newcommand{\revise}[1]{{\color{black} #1}}

\title{Hadamard's variational formula for simple eigenvalues}

\author{Takashi Suzuki
\thanks{Center for Mathematical Modeling and Data Science, Osaka University}
 \and \revise{Takuya} Tsuchiya\;${}^{*}$}

\date{}

\begin{document}

\maketitle

\begin{abstract}
We study Hadamard's variational formula for simple eigenvalues under dynamical and conformal deformations. Particularly, harmonic convexity of the first eigenvalue of the Laplacian under the mixed boundary condition is established for two-dimensional domain, which implies several new inequalities.       
\end{abstract}

\section{Introduction}\label{sec1}

The purpose of the present paper is to derive Hadamard's \revise{variational} formulae for simple eigenvalue of the Laplacian and is to derive several new inequalities concerning the first eigenvalue \revise{concerning} the mixed boundary condition. 

In the previous work \cite{st24} we studied Hadamard's variational formula for general domain deformation and extended the results \cite{gs52} on two or three-space dimensions under normal perturbations of the domain. There, we develop an abstract theory \revise{of} perturbation of self-adjoint operators, refining the argument \revise{in} \cite{rellich53}. 

To be precise, let $(X, \vert \ \vert)$ and $(V, \Vert \ \Vert)$ be a pair of Hilbert spaces over $\R$ with compact embedding $V\hookrightarrow X$. Henceforth, $C>0$ denotes a generic constant. 
Let $A_t:V\times V\rightarrow \R$ and $B_t:X\times X\rightarrow \R$ for $t\in I=(-\varepsilon_0, \varepsilon_0)$, $\varepsilon_0>0$, be symmetric bilinear forms satisfying 
\begin{equation}  
\vert A_t(u,v)\vert \leq C\Vert u\Vert\Vert v\Vert, \ A_t(u,u)\geq \delta \Vert u\Vert^2, \quad u,v\in V 
 \label{coercive1}
\end{equation} 
and 
\begin{equation} 
\vert B_t(u,v)\vert \leq C\vert u\vert\vert v\vert, \ B_t(u,u)\geq \delta\vert u\vert^2, \quad u, v\in X 
 \label{coercive2}
\end{equation} 
for some $\delta>0$.  We take the abstract eigenvalue problem  
\begin{equation} 
u\in V, \ B_t(u,u)=1, \quad A_t(u,v)=\lambda B_t(u,v), \ \forall v\in V,  
 \label{3}
\end{equation} 
which ensures a sequence of eigenvalues denoted by 
\[ 0<\lambda_1(t)\leq \lambda_2(t)\leq \cdots \rightarrow +\infty. \] 
The associated normalized eigenfunctions, 
\[ u_1(t), u_2(t), \cdots, \] 
furthermore, form a complete ortho-normal system in $X$, provided with the inner product induced by $B_t=B_t(\cdot, \cdot)$. Hence it holds that 
\[ B_t(u_i(t), u_j(t))=\delta_{ij}, \ A_t(u_j(t), v)=\lambda_jB_t(u_j(t), v) \] 
for any $v \in V$ and $i,j=1,2,\cdots$. 

The abstract theory developed in \cite{st24} is stated as follows. First, \revise{if 
\begin{eqnarray*} 
& & \lim_{h\rightarrow 0}\sup_{\Vert u\Vert, \Vert v\Vert\leq 1}\vert A_{t+h}(u,v)-A_t(u,v)\vert=0 \\ 
& & \lim_{h\rightarrow 0}\sup_{\vert u\vert, \vert v\vert\leq 1}\vert B_{t+h}(u,v)-B_t(u,v)\vert =0 
\end{eqnarray*}  
holds for fixed $t\in I$, we obtain  
\[ \lim_{h\rightarrow 0}\lambda_j(t+h)=\lambda_j(t) \] 
at this $t$, where $j=1,2,\cdots$ are arbitrary (Theorem 8 of \cite{st24}). Second, if there exist bilinear forms $\dot A_t:V\times V\rightarrow \R$ and $\dot B_t:X\times X\rightarrow \R$ for any $t\in I$} such that   
\begin{eqnarray} 
& & \vert \dot A_t(u,v)\vert  \leq C\Vert u\Vert\Vert v\Vert, \quad u,v\in V \nonumber\\ 
& & \vert \dot B_t(u,v)\vert  \leq C\vert u\vert\vert v\vert, \quad u, v\in X,  
 \label{d11}  
\end{eqnarray} 
and if it holds that  
\begin{eqnarray} 
& & \lim_{h\rightarrow 0}\frac{1}{h}\sup_{\Vert u\Vert, \Vert v\Vert\leq 1}
\left\vert \left(A_{t+h}-A_t-h\dot A_t\right)(u,v)\right\vert =0 \nonumber\\ 
& & \lim_{h\rightarrow 0}\frac{1}{h}\sup_{\vert u\vert, \vert v\vert\leq 1}
\left\vert \left(B_{t+h}-B_t-h\dot B_t\right)(u,v)\right\vert =0, 
 \label{abd}
\end{eqnarray} 
for fixed $t\in I$, we obtain the existence of \revise{the} unilateral derivatives 
\begin{equation} 
\dot \lambda_j^\pm(t)=\lim_{h\rightarrow \pm 0}\frac{1}{h}\{ \lambda_j(t+h)-\lambda_j(t)\} 
 \label{49} 
\end{equation}  
at this $t$, where $j=1,2,\cdots$ are arbitrary (Theorem 12 of \cite{st24}). 
\revise{If inequalities (\ref{abd}) are valid to any $t\in I$ and it also holds that 
\begin{eqnarray} 
& & \lim_{h\rightarrow 0}\sup_{\Vert u\Vert_V, \Vert v\Vert_V\leq 1}\left\vert \dot A_{t+h}(u,v)-\dot A_t(u,v)\right\vert=0 \nonumber\\ 
& & \lim_{h\rightarrow 0}\sup_{\vert u\vert_X, \vert v\vert_X\leq 1}\left \vert \dot B_{t+h}(u,v)-\dot B_t(u,v)\right\vert=0    
 \label{72x}
\end{eqnarray} 
for the specified $t$, furthermore, the above unilateral derivatives satisfy 
\[ \lim_{h\rightarrow \pm 0}\dot \lambda_j^\pm(t+h)=\dot \lambda_j^\pm (t) \]  
(Theorem 13 of \cite{st24}). } 

We assume, furthermore, 
\begin{equation}   
\lambda_{k-1}(t)<\lambda_k(t)\leq \cdots \leq \lambda_{k+m-1}(t)<\lambda_{k+m}(t), \quad \forall t\in I 
 \label{8}
\end{equation} 
for some $k,m=1,2,\cdots$, under the agreement of $\lambda_{0}(t)=-\infty$. Assume, also that (\ref{abd}) and (\ref{72x}) hold for any $t\in I$. Then, there exists a family of $C^1$ curves denoted by $\tilde C_i$, $k\leq i\leq k+m-1$, made by a rearrangement of 
\[ \{ C_j\mid k\leq j\leq k+m-1\} \] 
at most countably many times in $I$, where   
\[ C_j=\{ \lambda_j(t) \mid t\in I\}, \quad k\leq j\leq k+m-1 \] 
(Theorem 3, Theorem 14 of \cite{st24}). 

\revise{Third}, if we have the other bilinear forms $\ddot A_t:V\times V\rightarrow \R$ and $\ddot B_t:X\times X\rightarrow \R$ satisfying 
\begin{eqnarray} 
& & \vert \ddot A_t(u,v)\vert \leq C\Vert u\Vert\Vert v\Vert, \quad u,v\in V \nonumber\\ 
& & \vert \ddot B_t(u,v)\vert \leq C\vert u\vert\vert v\vert, \quad \ u,v\in X 
 \label{abd2+} 
\end{eqnarray} 
for any $t\in I$, and  
\begin{eqnarray} 
& & \lim_{h\rightarrow 0}\frac{1}{h^2}\sup_{\Vert u\Vert, \Vert v\Vert\leq 1}\left\vert \left(A_{t+h}-A_t-h\dot A_t-\frac{h^2}{2}\ddot A_t\right)(u,v)\right\vert=0 \nonumber\\ 
& & \lim_{h\rightarrow 0} \frac{1}{h^2}\sup_{\vert u\vert, \vert v\vert \leq 1}\left\vert \left( B_{t+h}-B_t-h\dot B_t-\frac{h^2}{2}\ddot B_t\right)(u,v)\right\vert =0 
 \label{abd2}
\end{eqnarray} 
for the fixed $t$, then there are  
\[ \ddot\lambda_j^\pm(t)=\lim_{h\rightarrow 0}\frac{2}{h^2}(\lambda_j(t+h)-\lambda_j(t)-h\dot \lambda_j^\pm(t)) \]  
for this $t$, where $j=1,2,\cdots$ are arbitrary (Remark 12 of \cite{st24}). These $\ddot \lambda_j^\pm(t)$, furthermore, satisfy  
\[ \lim_{h\rightarrow \pm 0}\ddot \lambda_j^\pm(t+h)=\ddot \lambda_j^\pm(t), \]  
if inequalities (\ref{abd2}) are valid to any $t\in I$ and it holds that 
\begin{eqnarray} 
& & \lim_{h\rightarrow 0}\sup_{\Vert u\Vert, \Vert v\Vert\leq 1}\left\vert \ddot A_{t+h}(u,v)-\ddot A_t(u,v)\right\vert=0 \nonumber\\ 
& & \lim_{h\rightarrow 0}\sup_{\vert u\vert, \vert v\vert\leq 1}\left\vert \ddot B_{t+h}(u,v)-\ddot B_t(u,v)\right\vert=0 
 \label{abdot}
\end{eqnarray} 
for the above specified $t$ (Theorem 23 of \cite{st24}).  \revise{Furthermore}, if (\ref{abd2}) and (\ref{abdot}) hold for any $t\in I$ and if it holds that (\ref{8}), then the above described family of curves $\tilde C_i$,  $k\leq i\leq k+m-1$, are $C^2$ (Theorem 24 of \cite{st24}). 

\revise{Finally, these} derivatives $\dot\lambda_j^\pm(t)$ and $\ddot \lambda_j^\pm(t)$ for $k\leq j\leq k+m-1$ are characterized as the eigenvalues of associated eigenvalue problems on the $m$-dimensional space, $\langle u_j(t) \mid k\leq j\leq k+m-1\rangle$ (Theorem 12 and Theorem 15 of \cite{st24}). Henceforth, we assume (\ref{coercive1}), (\ref{coercive2}), (\ref{d11}), (\ref{abd}), (\ref{72x}), (\ref{abd2+}), (\ref{abd2}), and (\ref{abdot}) for any $t\in I$. 

If $m=1$ \revise{holds} in (\ref{8}), for example, we obtain 
\begin{equation} 
\dot\lambda\equiv \dot\lambda_k^\pm(t)=(\dot A_t-\lambda \dot B_t)(\phi_t, \phi_t) 
 \label{dot}
\end{equation} 
and 
\[\ \ddot \lambda\equiv \ddot \lambda_k^\pm(t)=(\ddot A_t-\lambda \ddot B_t-2\dot \lambda B_t)(\phi_t, \phi_t)-2\revise{C_t^\lambda}(\gamma(\phi_t), \gamma(\phi_t)), \] 
where $\lambda=\lambda_k(t)$, $\phi_t=u_k(t)$, and $C_t^\lambda=A_t-\lambda B_t$. Here, $w=\gamma(u)\in V_1$ is defined for $u\in V$ by 
\[ C_t^\lambda(w,v)=-\dot C_t^{\lambda, \dot \lambda}(u,v), \quad \forall v\in V_1, \] 
where $V_1=PV$,  $P=I-R$, and $R:X\rightarrow V_0=\langle \phi_t\rangle$ is the orthogonal projection, and  
\[ \dot C_t^{\lambda, \dot \lambda}=\dot A_t-\lambda\dot B_t-\dot \lambda B_t. \] 
Thus we obtain the following theorem because the bilinear form $C_t^\lambda$ is non-negative definite on $V\times V$ if $\lambda=\lambda_1$. 

\begin{thm}
If $k=m=1$ \revise{holds} in (\ref{8}), \revise{there arises} that 
\begin{equation} 
\dot \lambda=\dot A-\lambda \dot B, \quad \ddot \lambda\leq \ddot A-\lambda \ddot B-2\dot\lambda\dot B, 
 \label{x7}
\end{equation} 
where 
\begin{equation} 
\dot A=\dot A_t(\phi_t,\phi_t), \ \ddot A=\ddot A_t(\phi_t, \phi_t), \ \dot B=\dot B_t(\phi_t, \phi_t), \ \ddot B=\ddot B_t(\phi_t, \phi_t) 
 \label{12a}
\end{equation} 
for $\phi_t=u_k(t)$. 
 \label{thm1}
\end{thm}

Since 
\[ \frac{d^2}{dt^2}\frac{1}{\lambda}=-\frac{\ddot \lambda \lambda^2-2\lambda \dot\lambda^2}{\lambda^2} \] 
Theorem \ref{thm1} implies the following result. 

\begin{thm} 
If $k=m=1$ and 
\begin{equation} 
2\dot A^2+\lambda^2\ddot B\geq \lambda(\ddot A+2\dot A\dot B), 
 \label{12}
\end{equation} 
it holds that 
\begin{equation} 
\frac{d^2}{dt^2}\frac{1}{\lambda}\geq 0. 
 \label{13}
\end{equation} 
 \label{thm2}
\end{thm} 

Harmonic convexity of the first eigenvalue, inequality (\ref{13}), was noticed by \cite{gs52} for the Dirichlet problem of \revise{Laplacian} under the conformal deformations of the domain in two space dimension. Here we calculate the values $\dot A$, $\ddot A$, $\dot B$, $\ddot B$ and the validity of (\ref{12}) under general setting of the deformations of the domain, and then turn to the dynamical and the conformal deformations. Taking preliminaries in $\S$\ref{sec2}, thus, we show the results on dynamical and conformal deformations in $\S$\ref{sec3} \revise{(Theorem 5 and Theorem 6)} and $\S$\ref{sec4} \revise{(Theorem 8 and Theorem 9)}, respectively. As applications, we show several new inequalities on the first eigenvalue of the two-dimensional problem.

\section{General deformations around $t=0$}\label{sec2} 

Let $\Omega$ be a bounded Lipschitz domain in $n$-dimensional Euclidean space $\R^n$ for $n \ge 2$. Suppose that its boundary $\partial\Omega$ is divided into two relatively open  \revise{disconnected}  
sets $\gamma_0$ and $\gamma_1$, satisfying 
\begin{equation}
   \overline{\gamma_0} \cup \overline{\gamma_1}= \partial\Omega,
    \quad
   \overline{\gamma_0} \cap \overline{\gamma_1} = \emptyset.
  \label{assum-gamma}
\end{equation}

We study the eigenvalue problem of the Laplacian with mixed boundary condition, 
\begin{equation}
  - \Delta u = \lambda u \ \mbox{in $\Omega$}, \quad
      u = 0 \ \mbox{on $\gamma_0$}, \quad
      \frac{\partial u}{\partial \nu} = 0 \ \mbox{on $\gamma_1$},
   \label{2}
\end{equation}
where
\[ \Delta=\sum_{i=1}^n \frac{\partial^2}{\partial x_i^2} \] 
and $\nu$ denotes the outer unit normal vector on $\partial\Omega$. This problem takes the weak form, finding $u$ satisfying  
\begin{equation} 
u\in V, \ B(u,u)=1, \quad A(u,v)=\lambda B(u,v), \ \forall v\in V 
 \label{3+}
\end{equation} 
defined for 
\begin{equation}  
A(u,v)=\int_\Omega \nabla u\cdot \nabla v \ dx, \quad u,v\in V 
 \label{auv}
\end{equation} 
and 
\begin{equation} 
B(u,v)=\int_\Omega uv \ dx,\quad u, v\in X,  
 \label{buv}
\end{equation}  
where  
\begin{equation} 
X=L^2(\Omega), \quad V=\{ v\in H^1(\Omega) \mid \left. v\right\vert_{\gamma_0}=0\}.   
 \label{v}
\end{equation} 
This $V$ is a closed subspace of $H^1(\Omega)$ under the norm 
\[ \Vert v\Vert=\left( \int_\Omega \vert \nabla v\vert^2+v^2 \ dx \right)^{1/2}. \]  
The above reduction of (\ref{2}) to \revise{(\ref{3+})}-(\ref{v}) is justified via the trace operator to the boundary since $\Omega$ is a bounded Lipschitz domain (Theorem 2 of \cite{st23}).  

To confirm the well-posedness of (\ref{3}), we note, first, that 
if $\gamma_0\neq\emptyset$, there is a coercivity of $A:V\times V\rightarrow \R$, which means the existence of $\delta>0$ such that  
\begin{equation} 
A(v,v)\geq \delta \Vert v\Vert^2, \quad \forall v\in V. 
 \label{4}
\end{equation} 
If $\gamma_0=\emptyset$ we replace $A$ by $A+B$, denoted by $\tilde A$. Then this $\tilde A:V\times V\rightarrow \R$ is coercive, and the eigenvalue problem 
\[ u\in V, \ B(u,u)=1, \quad \tilde A(u,v)=\tilde \lambda B(u,v), \ \forall v\in V, \] 
is equivalent to (\ref{3}) by $\tilde \lambda=\lambda+1$.  Hence  we can assume (\ref{coercive1}), using this  reduction if it is necessary. 

If the bounded domain $\Omega\subset {\bf R}^n$ is provided with the cone property, \revise{the inclusion} $H^1(\Omega) \hookrightarrow L^2(\Omega)$ is compact by Rellich-Kondrachov's theorem \cite{af05}. Thus there is a sequence of eigenvalues to (\ref{3}), denoted by 
\[ 0<\lambda_1\leq \lambda_2\leq \cdots \rightarrow +\infty. \] 
The associated eigenfunctions, $u_1, u_2, \cdots$, furthermore, form a complete ortho-normal system in $X$, provided with the inner product induced by $B=B(\cdot, \cdot)$: 
\[ B(u_i, u_j)=\delta_{ij}, \quad A(u_j, v)=\lambda_jB(u_j, v), \ \forall v \in V, \quad i,j=1,2,\cdots. \]  

The $j$-th eigenvalue of (\ref{3}) is given by the mini-max principle 
\begin{equation} 
\lambda_j=\min_{L_j}\max_{v\in L_j\setminus \{0\}}R[v]=\max_{W_j}\min_{v\in W_j\setminus \{0\}}R[v], 
 \label{minimax}
\end{equation} 
where 
\[ R[v]=\frac{A(v,v)}{B(v,v)} \] 
is the Rayleigh quotient, and $\{ L_j\}$ and $\{ W_j\}$ denote the families of all subspaces of $V$ with dimension and codimension $j$ and $j-1$, respectively.     

\revise{ The following well-known fact is valid without the smoothness of $\partial\Omega$. The proof given in Appendix A for completeness.} 
\begin{thm}
If $\Omega\subset{\bf R}^n$ is a bounded Lipschitz domain, the first eigenvalue $\lambda_1$ to (\ref{3}) formulated to $A(u,v)$ and $B(u,v)$ defined by (\ref{auv})-(\ref{v}) is simple.   
 \label{theorem3} 
\end{thm} 

Coming back to the Lipschitz bounded domain $\Omega$, we introduce its deformation as follows. Let  
\begin{equation} 
T_t : \Omega\rightarrow \Omega_t= T_t(\Omega), \quad  t\in (-\varepsilon_0, \varepsilon_0)    
 \label{13+}
\end{equation} 
be a family of bi-Lipschitz homeomorphisms. We assume that $T_tx$ is continuous in $t$ uniformly in $x\in \Omega$, and continue to use the following definition as in \cite{st24}.  

\begin{de}\label{def1} 
The family  $\{ T_t\}$ of bi-Lipschitz homeomorphisms is said to be $p$-differentiable in $t$ for $p\geq 1$, if $T_tx$ is  $p$-times differentiable in $t$ for any $x\in\Omega$ and the mappings 
\[ 
\frac{\partial^\ell}{\partial t^\ell}DT_t, \ \frac{\partial^\ell}{\partial t^\ell}(DT_t)^{-1}:\Omega \rightarrow M_n(\R), \quad 0\leq \ell \leq p  \nonumber\\  
\] 
are uniformly bounded in $(x,t)\in \Omega\times (-\varepsilon_0, \varepsilon_0)$, where $DT_t$ denotes the Jacobi matrix of $T_t:\Omega \rightarrow \Omega_t$ and $M_n(\R)$ stands for the set of real $n\times n$ matrices. This $\{ T_t\}$ is \revise{furthermore} said to be continuously $p$-differentiable in $t$ if it is $p$-differentiable and the mappings 
\[ 
t\in (-\varepsilon_0, \varepsilon_0)\mapsto \frac{\partial^\ell}{\partial t^\ell}DT_t, \ \frac{\partial^\ell}{\partial t^\ell}(DT_t)^{-1} \in L^\infty(\Omega \rightarrow M_n(\R)), \quad 0\leq \ell \leq p  \] 
are continuous.  
\end{de}

Putting  
\begin{equation} 
T_t(\gamma_i)=\gamma_{it}, \quad i=0,1, 
 \label{6}
\end{equation} 
in (\ref{2}), we introduce the other eigenvalue problem 
\begin{equation}
  - \Delta u= \lambda u \ \mbox{in $\Omega_t$}, \quad
      u = 0 \ \mbox{on $\gamma_{0t}$}, \quad
      \frac{\partial u}{\partial \nu} = 0 \ \mbox{on $\gamma_{1t}$},
   \label{7}
\end{equation}
which is reduced to finding 
\begin{equation} 
u\in V_t, \ \int_{\Omega_t}u^2 \ dx=1, \quad \int_{\Omega_t}\nabla u\cdot \nabla v \ dx=\lambda\int_{\Omega_t}uv \ dx, \ \forall v\in V_t 
 \label{88}
\end{equation} 
for 
\begin{equation} 
V_t=\{ v\in H^1(\Omega_t) \mid \left. v\right\vert_{\gamma_{0t}}=0\}.  
 \label{9}
\end{equation} 
Let $\lambda_j(t)$ be the $j$-th eigenvalue of the eigenvalue problem (\ref{7}). Then Lemma 7 of \cite{st24} ensures that the eigenvalue problem (\ref{88})-(\ref{9}) is reduced to 
\begin{equation} 
u\in V, \ B_t(u,u)=1, \quad A_t(u,v)=\lambda B_t(u,v), \ \forall v\in V 
 \label{wf}
\end{equation} 
by the transformation of variables $y=T_tx$, where 
\begin{eqnarray} 
& & B_t(u,v)=\int_\Omega uv a_t \ dx, \quad u, v\in X \nonumber\\ 
& & A_t(u,v)=\int_\Omega Q_t[\nabla u, \nabla v] a_t \ dx, \quad u,v\in V 
 \label{ab}
\end{eqnarray} 
for $V$ and $X$ defined by (\ref{v}), and 
\begin{equation} 
a_t=\det DT_t, \quad Q_t=(DT_t)^{-1}(DT_t)^{-1T}.  
 \label{cd}
\end{equation}

Recall that $\Omega\subset \R^n$ is a bounded Lipshitz domain, and let $T_t:\Omega\rightarrow \Omega_t=T_t\Omega$, $t\in I=(-\varepsilon_0, \varepsilon_0)$ be a family of twice continuously differentiable bi-Lispchitz transformations. Then the abstract theory described in $\S$\ref{sec1} is applicable with 
\[ \dot A_t(u,v)=\int_\Omega \frac{\partial}{\partial t}(Q_t[\nabla u, \nabla v]a_t)dx, \quad \dot B_t(u,v)=\int_\Omega uv \ \frac{\partial a_t}{\partial t} \ dx \] 
and 
\[ \ddot A_t(u,v)=\int_\Omega \frac{\partial^2}{\partial t^2}(Q_t[\nabla u, \nabla v]a_t)dx, \quad \ddot B_t(u,v)=\int_\Omega uv \ \frac{\partial^2 a_t}{\partial t^2} \ dx \] 

Henceforth, we write 
\[ G[\xi,\eta]=\xi^TG\eta, \quad \xi,\eta \in \R^n \] 
for the symmetric matrix $G\in M_n(\R)$, where ${\cdot}^T$ denotes the transpose \revise{of vectors or matrices}. The unit matrix is denoted by $E\in M_n(\R)$, and 
\[ G:F=\sum_{i,j=1}^ng_{ij}f_{ij} \] 
for $G=(g_{ij})$ and $F=(f_{ij})\in M_n({\bf R})$. Let, \revise{furthermore}, $I:\Omega\rightarrow \Omega$ be the identity mapping. 

\begin{thm} 
Assume   
\begin{equation} 
T_t=I+tS+\frac{t^2}{2}R+o(t^2), \quad t\rightarrow 0 
 \label{30a}
\end{equation} 
uniformly on $\Omega$, where $S,R:\Omega\rightarrow \R^n$ are Lipschitz continuous vector fields. Then, if $m=1$ in (\ref{8}) for $t=0$, it holds that  
\begin{eqnarray*} 
& & \dot A=\int_\Omega (-(DS^T+DS)+(\nabla\cdot S)E)[\nabla\varphi, \nabla\varphi] \\ 
& & \ddot A=\int_\Omega((2(DS)^2-DR)^T+(2(DS)^2-DR) \\ 
& &  \quad +2\revise{(DS)(DS)^T} -(\nabla\cdot S)(DS^T+DS) \\ 
& & \quad +(\nabla\cdot R+(\nabla\cdot S)^2-DS^T:DS)E)[\nabla\phi, \nabla \phi] \ dx \\ 
& & \dot B=\int_\Omega (\nabla\cdot S)\phi^2 \ dx  \\ 
& & \ddot B=\int_\Omega (\nabla\cdot R+(\nabla\cdot S)^2-DS^T:DS)\phi^2 \ dx, 
\end{eqnarray*} 
where $\dot A$, $\ddot A$, $\dot B$, and $\ddot B$ are defined by (\ref{12a}) with $t=0$, and $\phi=\phi_0$.   
 \label{theorem4} 
\end{thm} 

{\it Proof:} \ First, since 
\[ DT_t=E+tDS+\frac{t^2}{2}DR+o(t^2) \] 
uniformly on $\Omega$, it holds that 
\[ (DT_t)^{-1}=E-tDS^T+\frac{t^2}{2}(2(DS)^2-DR)+o(t^2) \] 
uniformly on $\Omega$, which implies 
\[ a_t=\det DT_t=1+t\revise{\nabla\cdot S}+\frac{t^2}{2}(\nabla\cdot R+(\nabla\cdot S)^2-DS^T:DS)+o(t^2) \] 
\revise{uniformly} on $\Omega$ (Lemma 5 and Lemma 6 of \cite{st23}). 

Second, we have 
\[ Q_t[\nabla u, \nabla v]=(DT_t^{-1}\nabla u, DT_t^{-1T}\nabla v), \quad u,v \in V, \] 
using the standard $L^2$ inner product $( \ , \ )$, and therefore, it follows that 
\begin{eqnarray*} 
& & \left. \frac{\partial Q_t}{\partial t}[\nabla u, \nabla v]\right\vert_{t=0} \\ 
& & \quad = \left. (\frac{\partial}{\partial t}DT_t^{-1T}\nabla u, DT_t^{-1T}\nabla v)+(DT_t^{-1T}\nabla u, \frac{\partial}{\partial t}DT_t^{-1T}Dv) \right\vert_{t=0} \\ 
& & \quad = (-DS^T\nabla u, \nabla v)+(\nabla u, -DS^T\nabla v)=-(DS^T+DS)[\nabla u, \nabla v] 
\end{eqnarray*} 
and 
\begin{eqnarray*} 
& & \left. \frac{\partial^2Q_t}{\partial t^2}[\nabla u, \nabla v]\right\vert_{t=0} = \left. (\frac{\partial^2}{\partial t^2}DT_t^{-1T}\nabla u, DT_t^{-1T}\nabla v)\right. \\ 
& & \quad \quad \left. +2(\frac{\partial}{\partial t}DT_t^{-1T}\nabla u, \frac{\partial}{\partial t}DT_t^{-1T}\nabla v)+(DT_t^{-1T}\nabla u, \frac{\partial^2}{\partial t^2}DT_t^{-1T}Dv) \right\vert_{t=0} \\ 
& & \quad = (2(DS)^2-DR)^T\nabla u, \nabla v)+2(DS^T\nabla u, DS^T\nabla v) \\ 
& & \qquad +(\nabla u, (2(DS)^2-DR)^T\nabla v) \\ 
& & \quad =((2(DS)^2-DR)^T+(2(DS)^2-DR)+2DSDS^T)[\nabla u, \nabla v]. 
\end{eqnarray*} 

By these equalities, finally, we obtain 
\begin{eqnarray*} 
& & \left.\frac{\partial}{\partial t}(Q_t[\nabla u, \nabla v]a_t)\right\vert_{t=0}=\left. \frac{\partial Q_t}{\partial t}[\nabla u, \nabla v]a_t+Q_t[\nabla u, \nabla v]\frac{\partial a_t}{\partial t}\right\vert_{t=0} \\ 
& & \quad =-(DS^T+DS)[\nabla u, \nabla v]+(\nabla u, \nabla v)\nabla\cdot S 
\end{eqnarray*} 
and 
\begin{eqnarray*} 
& & \left.\frac{\partial^2}{\partial t^2}(Q_t[\nabla u, \nabla v]a_t)\right\vert_{t=0}=\left. \frac{\partial^2Q_t}{\partial t^2}[\nabla u, \nabla v]a_t \right. \\ 
& & \qquad \left. +2\frac{\partial Q_t}{\partial t}[\nabla u, \nabla v]\frac{\partial a_t}{\partial t}+Q_t[\nabla u, \nabla v]\frac{\partial^2 a_t}{\partial t^2}\right\vert_{t=0} \\ 
& & \quad =-((2(DS)^2-DR)+(2(DS)^2-DR)+2DSDS^T)[\nabla u, \nabla v] \\ 
& & \qquad -2(DS^T+DS)[\nabla u, \nabla v](\nabla \cdot S) \\ 
& & \qquad +(\nabla u, \nabla v)(\nabla\cdot R+(\nabla \cdot S)^2-DS^T:DS),  
\end{eqnarray*} 
and hence the conclusion. 

\section{Dynamical deformations}\label{sec3}

We continue to suppose that $\Omega \subset {\bf R}^n$ is a bounded Lipschitz domain. Here we study the dynamical deformation of domains introduced by \cite{st16}.  

To this end, we take a Lipschitz continuous vector field defined on a neighbourhood of $\overline{\Omega}$, denoted by $v=v(x)$. Then the transformation $T_t:\Omega\rightarrow \Omega_t=T_t(\Omega)$ is made by $T_tx=X(t)$, where $X=X(t)$, $\vert t\vert \ll 1$, is the solution to 
\[ \frac{dX}{dt}=v(X), \quad \left. X\right\vert_{t=0}=x. \] 
Then we have \revise{the group property}, 
\[ T_t\circ T_s=T_{t+s}, \quad \vert t\vert, \vert s\vert\ll 1, \] 
and therefore, the formulae derived for $t=0$ are shifted to \revise{$t=s$.}  

If $v$ is a $C^{1,1}$ vector field, furthermore, this $T_t:\Omega\rightarrow \Omega_t$ is twice continuously differentiable, and it holds that  (\ref{30a}) with 
\begin{equation}  
S=v, \quad R=(v\cdot \nabla)v 
 \label{x1} 
\end{equation} 
because of 
\[ \frac{d^2X}{dt^2}=[(v\cdot \nabla)]v(X).  \] 

Then we obtain the following lemma. 

\begin{lem} 
Under the above assumption, it holds that 
\[ \nabla\cdot R=DS^T:DS+(v\cdot \nabla)(\nabla\cdot v). \]
 \label{lemma1}
\end{lem} 
{\it Proof:} \ Writing $v=(v^j)$, we obtain 
\begin{eqnarray*} 
\nabla \cdot R & = & \nabla\cdot ((v\cdot \nabla)v) =\sum_{i,j}\frac{\partial}{\partial x_i}(v^j\frac{\partial v^i}{\partial x_j}) \\ 
& = & \sum_{i,j}(\frac{\partial v^j}{\partial x_i}\frac{\partial v^i}{\partial x_j}+v^j\frac{\partial}{\partial x_j}\frac{\partial v^i}{\partial x_i}) = DS^T:DS+(v\cdot \nabla)(\nabla\cdot v)   
\end{eqnarray*} 
and the proof is complete. \hfill $\Box$ 

\vspace{2mm} 

Here we take two categories that the vector fields are solenoidal and gradient. In the first category, we assume $\nabla\cdot v=0$ everywhere. Then it holds that \revise{ 
\begin{equation} 
\vert \Omega_t\vert=\vert \Omega\vert, \quad \vert t\vert \ll 1. 
 \label{35+}
\end{equation} 
}

\begin{thm}
If $\nabla\cdot v=0$, it holds that $\dot B=\ddot B=0$ and 
\begin{eqnarray} 
& & \dot A=-\int_\Omega (DS^T+DS)[\nabla\phi, \nabla\phi] \ dx \nonumber\\ 
& & \ddot A=\int_\Omega (2(DS)^{2}-DR)^T+(2(DS)^2-DR) \nonumber\\ 
& & \qquad +2\revise{(DS)(DS)^T})[\nabla\phi, \nabla\phi] \ dx.  
 \label{32}  
\end{eqnarray}  
in Theorem \ref{theorem4}. 
 \label{theorem5}
\end{thm} 
{\it Proof:} \ We recall (\ref{x1}). By the assumption we obtain $\nabla\cdot S=0$ and 
\[ \nabla\cdot R=DS^T:DS \] 
by Lemma \ref{lemma1}, which implies $\dot B=\ddot B=0$, and (\ref{32}) by Theorem \ref{theorem4}. \hfill $\Box$ 

\begin{rem} 
For the moment we write 
\[ a_i=\frac{\partial a}{\partial x_i}, \quad a_{ij}=\frac{\partial^2a}{\partial x_i\partial x_j} \] 
in short. Then, by the proof of Lemma \ref{lemma1} we obtain 
\[ DR=\revise{((v^iv^j_i)_k)_{jk}=\revise{(DS)^2}+(v^iv^j_{ik})_{jk}} \] 
and 
\[ DR^T=\revise{(DS)^{2T}+(v^iv^k_{ij})_{jk}} \] 
similarly. Hence it holds that 
\[ \ddot A=\int_\Omega (\revise{(DS)^{2T}+(DS)^2+2(DS)(DS)^T}-G)[\nabla\phi, \nabla\phi] \ dx \] 
for 
\[ G=\revise{(v^i(v^j_{ik}+v^k_{ij}))_{jk}}. \] 
\end{rem} 

In the second category that the vector field is a gradient of a scalar field we obtain the following theorem. 

\begin{thm} 
If $v=\nabla \mu$ for the $C^{1,1}$ scalar field $\mu$, it holds that 
\begin{eqnarray} 
& & \dot A=\int_\Omega (-2H+(\mbox{tr} \ H)E)[\nabla \phi, \nabla \phi] \ dx \nonumber\\ 
& & \ddot A=\int_\Omega (6H^2-K-4(\mbox{tr} \ H)H \nonumber\\ 
& & \quad +((\mbox{tr} \ H)^2+ 
\frac{1}{2} \mbox{tr} \ K -H:H)E)[\nabla \phi, \nabla \phi] \ dx
 \label{x4}
\end{eqnarray} 
and 
\begin{equation} 
\dot B=\int_\Omega (\mbox{tr} \ H)\phi^2 \ dx, \quad \ddot B=\int_\Omega ((\mbox{tr} \ H)^2+\frac{1}{2}\mbox{tr} \ K-H:H)\phi^2 \ dx, 
 \label{x3}
\end{equation}  
where 
\[ H=\nabla^2 \mu, \quad K=\nabla^2(\vert \nabla\mu\vert^2) \] 
and 
\[ \vert H\vert^2=H:H=\sum_{ij}h_{ij}^2 \quad \mbox{for $H=(h_{ij})$}. \] 
 \label{theorem6} 
\end{thm} 

{\it Proof:} \ In this case, we obtain $S=\nabla \mu$ and 
\begin{equation} 
R=((\nabla \mu)\cdot \nabla)\nabla \mu=\left(\sum_k\mu_k\mu_{ik}\right)_i=\left( (\frac{1}{2}\sum_k\mu_k^2)_i\right)_i=\frac{1}{2}\nabla \vert \nabla \mu\vert^2. 
 \label{33}
\end{equation} 
Then it holds that 
\[ \dot B=\int_\Omega \phi^2\nabla\cdot S \ dx=\int_\Omega \phi^2\Delta \mu \ dx=\int_\Omega (\mbox{tr} \ H)\phi^2 \ dx. \] 
It also holds that 
\begin{eqnarray} 
(v\cdot \nabla)(\nabla\cdot v) & = & ((\nabla \mu)\cdot \nabla)(\Delta \mu)=\nabla \mu\cdot \Delta \nabla \mu \nonumber\\ 
& = & \frac{1}{2}\Delta (\vert \nabla \mu\vert^2)-\vert \nabla^2\mu \vert^2, 
 \label{x2}  
\end{eqnarray} 
because 
\[ \Delta \vert a\vert^2=2\Delta a\cdot a+\sum_i\frac{\partial a}{\partial x_i}\cdot \frac{\partial a}{\partial x_i} \] 
is valid to the vector field $a$. 

Since 
\[ \ddot B=\int_\Omega ((\nabla\cdot S)^2+(v\cdot\nabla)(\nabla\cdot v))\phi^2 \ dx \] 
is valid by Lemma \ref{lemma1}, and hence  
\begin{eqnarray*} 
\ddot B & = & \int_\Omega( (\Delta \mu)^2+\frac{1}{2}\Delta(\vert \nabla \mu\vert^2)-\vert \nabla^2\mu \vert^2)\phi^2 \ dx \\ 
& = & \int_\Omega (\mbox{tr} \ H)^2+ \frac{1}{2}\mbox{tr} \ K-\vert H\vert^2)\phi^2 \ dx. 
\end{eqnarray*} 

Next, we obtain 
\[ DS^T+DS=2H \] 
and hence 
\begin{eqnarray*} 
\dot A & = & \int_\Omega (-2H+(\Delta \mu)E)[\nabla \phi, \nabla \phi] \ dx \\ 
& = & \int_\Omega (-2H+(\mbox{tr} \ H)E)[\nabla\phi, \nabla\phi] \ dx. 
\end{eqnarray*} 

Finally, we divide $\ddot A$ as in 
\[ \ddot A=I+II+III \] 
for 
\begin{eqnarray*} 
& & I=\int_\Omega (2((DS)^{2T}+(DS)^2+DS^TDS))-(DR^T+DR))[\nabla\phi, \nabla\phi] \ dx \\ 
& & II=\int_\Omega -2(DS^T+DS)[\nabla\phi, \nabla\phi]\nabla\cdot S \ dx \\ 
& & III=\int_\Omega (\nabla\cdot R+(\nabla\cdot S)^2-DS^T:DS)\vert \nabla\phi\vert^2 \ dx. 
\end{eqnarray*} 
Here we have 
\[ (DS)^{2T}+(DS)^2+DS^TDS=3H^2. \] 
Equality (\ref{33}) now implies 
\[ R^i_j+R^j_i=(\vert \nabla \mu\vert^2)_{ij} \] 
for $R=(R^i)$, and hence 
\[ \nabla R^T+\nabla R=\nabla^2(\vert \nabla\mu\vert^2)=K. \] 
Thus we obtain 
\[ I=\int_\Omega (6H^2-K)[\nabla\phi, \nabla\phi] \ dx. \] 
It holds also that 
\[ II=\int_\Omega -\revise{4}H[\nabla\phi, \nabla\phi](\Delta \mu) \ dx =\int_\Omega -\revise{4}(\mbox{tr} \ H)H[\nabla\phi, \nabla\phi] \ dx. \] 
Finally, we have 
\begin{eqnarray*} 
III & = & \int_\Omega ((\nabla\cdot S)^2+(v\cdot \nabla)v)\vert \nabla\phi\vert^2 \ dx \\ 
& = & \int_\Omega ((\Delta \mu)^2+\frac{1}{2}\Delta(\vert \nabla\mu\vert^2)-\vert \nabla^2\mu\vert^2)\vert \nabla\phi\vert^2 \ dx \\ 
& = & \int_\Omega ((\mbox{tr} \ H)^2+\frac{1}{2}\mbox{tr} \ K-\vert H\vert^2)\vert \nabla \phi\vert^2 \ dx. 
\end{eqnarray*}  
We thus end up with the result by combining these equalities. \hfill $\Box$ 

\begin{rem} 
Theorem \ref{theorem5} and Theorem \ref{theorem6} are consistent if $v=\nabla \mu$ with 
\[ \Delta \mu=0.  \] 
In this case it holds that 
\begin{equation}  
\mbox{tr} \ H=0 
 \label{x8}
\end{equation} 
and 
\begin{equation} 
\frac{1}{2}\mbox{tr} \ K-\vert H\vert^2=0 
 \label{x9}
\end{equation} 
by (\ref{x2}): 
\[ 0=(v\cdot\nabla )(\nabla\cdot v)=\frac{1}{2}\Delta (\vert \nabla \mu\vert^2)-\vert \nabla^2\mu\vert^2. \] 
Then we can reproduce the result of Theorem \ref{theorem5} from (\ref{x3}), as in 
\begin{equation} 
\dot B=0, \quad \ddot B=0. 
 \label{x5}
\end{equation}  
We have, \revise{furthermore}, 
\begin{equation} 
\dot A=\int_\Omega -2H\vert \nabla \phi\vert^2 \ dx, \quad \ddot A=\int_\Omega (6H^2-K)\vert \nabla\phi\vert^2 \ dx 
 \label{x6}
\end{equation}  
by (\ref{x4}). By (\ref{x7}), (\ref{x5}) and (\ref{x6}), first, 
\revise{
\[ \dot \lambda=\dot A-\lambda \dot B=\dot A \] 
} 
cannot have a definite sign because of (\ref{x8}).   Second, if 
\begin{equation} 
6H^2-K\leq 0 
 \label{xx7}
\end{equation} 
we obtain \revise{ $\ddot A\leq 0$  by (\ref{x7}), which implies 
\[ \ddot \lambda\leq \ddot A-\lambda \ddot B-2\dot \lambda \dot B=\ddot A\leq 0 \] 
regardless of the form of $\phi$.} For inequality (\ref{xx7}) to hold,  however, it is necessary that  
\[ \mbox{tr} \ (6H^2-K)\geq 0, \] 
or 
\[ 3\mbox{tr} \ (H^2)\leq \vert H\vert^2 \] 
by (\ref{x9}). This inequality implies $H=0$ because 
\[ \mbox{tr} \ (H^2)=\sum_{ij}h_{ij}^2=\vert H\vert^2 \] 
for $H=(h_{ij})$ with $h_{ji}=h_{ij}$. In other words, several possibilities can arise to the monotonicity or \revise{convexity} of the first eigenvalue under the dynamical perturbation\revise{, provided with the property (\ref{35+}).} 
 
\end{rem} 

\section{Conformal deformations}\label{sec4}

Here we sayt that $T_t:\Omega\rightarrow \Omega_t$ is conformal if  
\[ (DT_t)^{-1}(DT_t)^{-1T}=\frac{1}{\alpha_t}E \] 
holds with a scalar field $\alpha_t>0$. It follows that $\alpha_tE=(DT_t)(DT_t)^T$, and therefore, the matrix $F_t=\alpha_t^{-1/2}DT_t$ is orthogonal: $F_tF_t^T=E$. 

Then we assume, furthermore, that 
\[ \det F_t=1 \] 
for simplicity. It follows that 
\begin{equation} 
a_t=\det DT_t=\det (\alpha_t^{1/2}F_t)=\alpha_t^{n/2}, 
 \label{jacobian}
\end{equation} 
and hence $\alpha_t=a_t^{2/n}$, which results in 
\[ (DT_t)^{-1}(DT_t)^{-1T}=a_t^{-2/n}E, \]  
and hence 
\[ A_t(u,v)=\int_\Omega (\nabla u\cdot \nabla v) a_t^{1-2/n}dx, \quad u,v\in V \] 
and 
\[ B_t(u,v)=\int_\Omega uva_t \ dx. \] 

Then we can show the following theorem. 

\begin{thm} \label{thm7}
Let $n=2$, $k=m=1$, and $T_t:\Omega\rightarrow \Omega_t$ be conformal at $t$. Assume, furthermore, that 
\[ \frac{\partial^2a_t}{\partial t^2}\geq 0 \quad \mbox{in $\Omega$}. \] 
Then it holds that 
\begin{equation} 
\frac{d^2}{dt^2}\frac{1}{\lambda_1(t)}\geq 0. 
 \label{hconcave}
\end{equation} 
 \label{theorem7}
\end{thm} 

{\it Proof:} \ Since $n=2$ we obtain 
\begin{equation} 
A_t(u,v)=\int_\Omega \nabla u\cdot \nabla v \ dx, \quad B_t(u,v)=\int_\Omega uva_t \ dx. 
 \label{abn2}
\end{equation} 
Then it holds that $\dot A=\ddot A=0$ and 
\[ \ddot B=\int_\Omega \phi^2\frac{\partial^2a_t}{\partial t^2} \ dx \] 
in Theorem \ref{thm1} for $\phi=\phi_t$. Then the result follows from Theorem \ref{thm2}. \hfill $\Box$ 

\vspace{2mm} 

A consequence of Theorem \ref{theorem7} is the following isometric inequality. 

\begin{thm} \label{thm8}
Let $\Omega \subset {\bf R}^2$ be a simply-connected \revise{bounded Lipschitz} domain and $f=f(z):D\rightarrow \Omega$ be a univalent \revise{bi-Lipschitz homeomorphism}, where $D=\{ z\in \mathbb{C} \mid \vert z\vert <1\}$. Let, furthermore, 
\[ f(z)=\sum_{n=0}^\infty a_nz^n, \quad a_1\in \R. \] 
Then it holds that 
\begin{equation} 
\lambda_1(D)\geq \lambda_1(\Omega)\left(2\int_D\phi_1^2\mbox{Re}f'(z) \ dx-1\right), 
 \label{34}
\end{equation}  
where $\lambda_1(\Omega)$ and $\lambda_1(D)$ are the first eigenvalue of (\ref{2}): 
\begin{equation} 
-\Delta u=\lambda u \ \mbox{in $\Omega$}, \quad u=0 \ \mbox{on $\gamma_0$}, \quad \frac{\partial u}{\partial \nu}=0 \ \mbox{on $\gamma_1$} 
 \label{evomega}
\end{equation} 
and that for $\Omega=D$, 
\begin{equation} 
-\Delta u=\lambda u \ \mbox{in $D$}, \quad u=0 \ \mbox{in $\gamma_0'$}, \quad \frac{\partial u}{\partial \nu}=0 \ \mbox{on $\gamma_1'$} 
 \label{36}
\end{equation} 
respectively with $\gamma_i=f(\gamma_i')$, $i=1,2$, and furthermore, $\phi_1=\phi_1(x)$ is the first eigenfunction of (\ref{36}) such that 
\[ \int_D \phi_1^2 \ dx=1. \]   
 \label{theorem8}
\end{thm} 

\begin{rem} 
\revise{Either $\gamma_0=\emptyset$ or $\gamma_1=\emptyset$ occurs if $\Omega\subset \R^2$ is simply-connected and  (\ref{assum-gamma}), and in} the former case it follows that $\lambda_1(\Omega)=0$ and the above theorem is trivial. This assumption (\ref{assum-gamma}), however, is used in \cite{st24} just to formulate  (\ref{2}) as in (\ref{3}) with (\ref{auv})-(\ref{buv}), and (\ref{v}). \revise{Hence we can exclude this assumption if we formulate (\ref{evomega}) as in (\ref{3+}) with (\ref{auv})-(\ref{buv}) for $V$ standing for the closure of $V_0$ in $H^1(\Omega)$,  where $v\in V_0$ if and only if there is $\tilde v$, a smooth extension of $v$ in a neighborhood of $\overline{\Omega}$, such that 
\[ \left. \tilde v\right\vert_{\gamma_0}=0. \]  
Then we re-formulate (\ref{36}) as in (\ref{wf}) with (\ref{ab}) for $t=1$, using the bi-Lipshichits homeomorphism $T_1=f^{-1}:\Omega\rightarrow D$. Under this agreement of (\ref{evomega}) and (\ref{36}), we exclude the assumption (\ref{assum-gamma}) in the above theorem. }  
\end{rem} 

{\it Proof of Theorem \ref{theorem8}:} \ Let 
\begin{equation} 
g_t(z)=(1-t)z+tf(z), \quad 0\leq t\leq 1 
 \label{35}
\end{equation} 
be the conformal mapping on $D$, and define $A_t$ and $B_t$ by (\ref{abn2}) for (\ref{jacobian}) with $n=2$. We  obtain $a_t=\vert g'_t(z)\vert^2$ in (\ref{jacobian}) and hence it follows that 
\[ \vert g_t'(z)\vert^2=\vert 1+t(f'(z)-1)\vert^2=1+2t \ \mbox{Re} \ (f'(z)-1)+t^2\vert f'(z)-1\vert^2 \] 
and hence 
\begin{eqnarray} 
& & \frac{\partial}{\partial t}\vert g_t'(z)\vert^2=2\mbox{Re}(f'(z)-1)+2t\vert f'(z)-1\vert^2 \nonumber\\ 
& & \frac{\partial^2}{\partial t^2}\vert g_t'(z)\vert^2=2\vert f'(z)-1\vert^2\geq 0.  
 \label{39}
\end{eqnarray} 
Thus we obtain 
\[ \frac{d^2}{dt^2}\frac{1}{\lambda_1(t)}\geq 0, \quad 0\leq t\leq 1 \] 
for $\lambda_1(t)$ defined by (\ref{minimax}) with $j=1$ and $A=A_t$ and $B=B_t$, which implies 
\begin{equation} 
\frac{1}{\lambda_1(\Omega)}\geq \frac{1}{\lambda_1(D)}+\left. \left(\frac{1}{\lambda_1(t)}\right)'\right\vert_{t=0}, 
 \label{41}
\end{equation} 
or 
\[ \frac{1}{\lambda_1(\Omega)}\geq \frac{1}{\lambda_1(D)} +\frac{1}{\lambda_1(D)}\int_D\varphi_1^2\left. \frac{\partial}{\partial t}\vert g_t'(z)\vert^2 \right\vert_{t=0} \ dx \] 
by 
\[ \left( \frac{1}{\lambda_1}\right)'=-\frac{\lambda_1'}{\lambda_1^2} \] 
and (\ref{dot}) with (\ref{abn2}). Then we obtain the result by (\ref{39}). \hfill $\Box$ 

\begin{rem} 
If $\gamma_1=\emptyset$, inequality (\ref{34}) is reduced to 
\begin{equation} 
\lambda_1(D)\geq \lambda_1(\Omega)(2\mbox{Re} \ a_1-1),
 \label{37}
\end{equation} 
Since $f:D\rightarrow \Omega$ is univalent, it holds that 
\[ \vert \Omega\vert=\int_D\vert f'(z)\vert^2 dx=\pi\sum_{n=1}^\infty n\vert a_n\vert^2, \] 
and therefore, $\vert\Omega\vert=\vert D\vert$ if and only if 
\[ \sum_{n=1}^\infty n\vert a_n\vert^2=1. \] 
This equality implies $\vert a_1\vert\leq 1$ \revise{and in particular, 
\[ 2\mbox{Re} \ a_1-1\leq 1 \] 
in (\ref{37}). } 
 \label{remark3}
\end{rem} 
 
Inequality (\ref{37}) for $\gamma_1=\emptyset$ is proven in Appendix B. We conclude this section with an analogous result to (\ref{34}). 

\begin{thm} 
Under the assumption of the previous theorem it holds that 
\[ \lambda_1(\Omega)\geq \lambda_1(D)\left( 2 \int_D\tilde \phi_1^2\mbox{Re}f'(z) \ dx-1\right), \] 
where $\tilde \phi_1=\hat\phi_1\circ f$ for the first eigenfunction $\hat \phi_1=\hat\phi_1(x)$ to (\ref{evomega}) such that 
\[ \int_\Omega \hat \phi_1^2 \ dx=1. \]  
 \label{theorem9} 
\end{thm}  

{\it Proof:} \ Similarly to (\ref{41}), we obtain 
\begin{eqnarray*} 
\frac{1}{\lambda_1(D)} & \geq & \frac{1}{\lambda_1(\Omega)}-\left.\left( \frac{1}{\lambda_1(t)}\right)'\right\vert_{t=1}=\frac{1}{\lambda_1(\Omega)}\left( 1+\frac{\lambda_1'(1)}{\lambda_1(\Omega)^2}\right) \\  
& = & \frac{1}{\lambda_1(\Omega)}\left( 1-\int_D\tilde \phi_1^2\left.\frac{\partial}{\partial t}\vert g_1'(z)\vert^2\right\vert_{t=1} \ dx\right),  
\end{eqnarray*} 
which implies the result by 
\[ \int_D\tilde \phi_1^2\vert f'(z)\vert^2 \ dx=\int_\Omega \hat \phi_1^2 \ dx=1 \] 
and 
\[ \left. \frac{\partial}{\partial t}\vert g'(z)\vert^2\right\vert_{t=1}=2\vert f'(z)\vert^2-2\mbox{Re}f'(z). 
\qquad \qquad \Box 
\]

\appendix

\section{Proof of Theorem \ref{theorem3}} 

\revise{Although using the Harnack inequality is standard (c.f. Theorem 8.38 of \cite{gt83}), here we follow the argument by \cite{ot88}. In fact, by }
(\ref{minimax}) it holds that 
\[ \lambda_1=\min_{v\in V\setminus \{0\}}R[v]. \] 
Since $V\hookrightarrow X$ is compact, there is $\phi_1\in V$ in $B(\phi_1, \phi_1)=1$ which attains $\lambda_1$. Then it holds that 
\begin{equation} 
A(\phi_1,v)=\lambda_1B(\phi_1, v), \quad \forall v\in V. 
 \label{a1}
\end{equation} 
We have 
\begin{equation} 
v\in H^1(\Omega) \ \Rightarrow \ [\max\{v,0\}]_{x_i}=
\left\{ \begin{array}{ll} 
v_{x_i}, & \{ x\in \Omega \mid v(x)>0\} \\ 
0,  & \{ x\in \Omega \mid v(x)\leq 0\} \end{array}\right. 
 \label{a2}
\end{equation} 
for $1\leq i\leq n$ (Definition 6.7 and Theorem A.1 of \cite{ks80}), and therefore, 
\[ v\in H^1(\Omega) \ \Rightarrow \ \vert v\vert\in H^1(\Omega) \] 
We thus obtain  
\begin{equation} 
v\in V \ \Rightarrow \ \vert v\vert \in V, \quad A(\vert v\vert, \vert v\vert)=A(v,v) 
 \label{a3}
\end{equation} 
because $\Omega$ is a Lipschitz domain. Hence we may assume $\phi_1\geq 0$ in $\Omega$. 

Equality (\ref{a1}) implies 
\[ -\Delta \phi_1=\lambda_1\phi_1 \quad \mbox{in $\Omega$} \] 
in the sense of distributions, and therefore, we obtain $\phi_1\in C^2(\Omega)$ by the elliptic regularity. Then the strong maximum principle implies $\phi_1>0$ in $\Omega$. 

Now we use the following lemma derived from (\ref{a2}). 

\begin{lem}\label{lem2} 
Given $v,w\in V$, let 
\[ m=\min\{ v,w\}, \quad M=\max\{ v,w\}. \] 
Then it follows that $m,M\in V$ and 
\begin{equation} 
A(m)+A(M)=A(v)+A(w), \quad B(m)+B(M)=B(v)+B(w). 
 \label{a4} 
\end{equation} 
\end{lem} 

{\it Proof:} \  We obtain $m,M\in H^1(\Omega)$ by (\ref{a2}) and hence $m,M\in V$ because $\Omega$ is a Lipschitz domain.  It is obvious that 
\begin{eqnarray*} 
B(m)+B(M) & = & \left(\int_{v>w}+\int_{v=w}+\int_{v<w}\right) \left(m^2+M^2\right) \ dx \\ 
& = & \left(\int_{v>w}+\int_{v=w}+\int_{v<w}\right) \left( v^2+w^2\right) \ dx \\ 
& = & B(v)+B(w),  
\end{eqnarray*} 
while 
\begin{eqnarray*} 
& & A(m)+A(M) = \left(\int_{v>w}+\int_{v=w}+\int_{v<w}\right)\left( \vert \nabla m\vert^2+\vert \nabla M\vert^2\right) \ dx \\ 
& & \quad = \left( \int_{v>w}+\int_{v=w}+\int_{v<w}\right)\left( \vert \nabla v\vert^2+\vert\nabla w\vert^2\right) \ dx \\ 
& & \quad = A(v)+A(w)  
\end{eqnarray*} 
follows from (\ref{a2}). \hfill $\Box$ 

\vspace{2mm} 

We are ready to give the following proof. 

\vspace{2mm} 

{\it Proof of Theorem \ref{theorem3}:} \ It sufficies to show 
\[ \nabla(z/\phi_1)=0 \quad \mbox{in $\Omega$}. \] 
for any first eigenfunction $z\in V$. To this end, we take $x_0\in \Omega$ and put 
\[ m=\min\{ z, t_0\phi_1\}, \quad M=\max\{ z, t_0\phi_1\} \] 
for $t_0=z(x_0)/\phi_1(x_0)$. The desired equality $\nabla (z/\phi_1)(x_0)=0$ is thus reduced to 
\[ \phi_1(x_0)\nabla z(x_0)=z(x_0)\nabla \phi_1(x_0), \] 
or, 
\begin{equation} 
\nabla z(x_0)=t_0\nabla \phi_1(x_0). 
 \label{a5}
\end{equation} 

Letting $C=A-\lambda_1B$, we obtain 
\[ C(m)+C(M)=C(z)+C(t_0\phi_1)=0 \] 
by Lemma \ref{lem2}, while $C$ is non-negative definite on $V\times V$. Hence it holds that $C(m)=C(M)=0$, and hence $m, M\in V$ are the other first eigenfunctions. We thus obtain $M\in C^2(\Omega)$, in particular. 

Given $e=(0,\cdots, \overbrace{1}^i, \cdots, 0)^T$, we obtain 
\[ z(x_0+he)-z(x_0)\leq M(x_0+he)-M(x_0), \quad \vert h\vert \ll 1 \] 
by $t_0\phi_1(x_0)=z(x_0)=M(x_0)$. Dividing both sides by $h$ and making $h\rightarrow \pm 0$, we thus obtain 
\[ z(x_0)_{x_i}=M(x_0)_{x_i}, \quad 1\leq i\leq n, \] 
and hence $\nabla M(x_0)=\nabla z(x_0)$. It holds that $\nabla
M(x_0)=t_0\nabla \phi_1(x_0)$, and hence (\ref{a5}). \hfill $\Box$

\section{Proof of (\ref{37}) for $\gamma_1=\emptyset$}

If $\gamma_1=0$, we have $\phi=\phi(r)$, and the result is a direct consequence of the following theorem. 

\begin{lem} 
If $\varphi=\varphi(r)$ for $r=\vert z\vert$ and $h=h(z)$ is holomorphic in $D$, it holds that 
\[ \int_D\varphi(r)^2h(z) \ dx= h(0)\int_D\phi(r)^2 \ dx. \] 
 \label{lem3}
\end{lem} 

{\it Proof of Lemma \ref{lem3}:} \  
Writing $z=re^{\imath \theta}$, we get  
\[ \int_D\varphi(r)^2h(z) \ dx=\int_0^1\varphi(r)^2rdr \cdot \int_0^{2\pi}h(re^{\imath\theta})d\theta. \] 
Since $dz=\imath re^{\imath\theta}d\theta=izd\theta$ on $\vert z\vert=1$, it holds that 
\[ \int_0^{2\pi}h(re^{\imath\theta})d\theta=\int_{\vert z\vert=r} \frac{h(z)}{\imath z} \ dz=2\pi h(0) \] 
because $h=h(z)$ is homeomorphic. Then we obatin 
\[ \int_D\phi(r)^2h(z) \ dx=2\pi h(0)\int_0^1\varphi(r)^2r \ dr = h(0)\int_D\varphi(r)^2 \ dx. \quad \Box \] 

\section{Numerical experiments}
In this section, we present numerical experiments to confirm the
theoretical results obtained.

Let $n = 2$ and $D \subset \R^2 \cong \C$ be the unit disk. 
Following \eqref{35},
we define $g_t: D \to \C$ by $g_t(z) = (1 - t)z + t \cos z$ for
$z \in \C$.  In this case, the images of the transformation are as
follows.

\begin{figure}[htbp] \label{fig1}
\includegraphics[width=2.0cm]{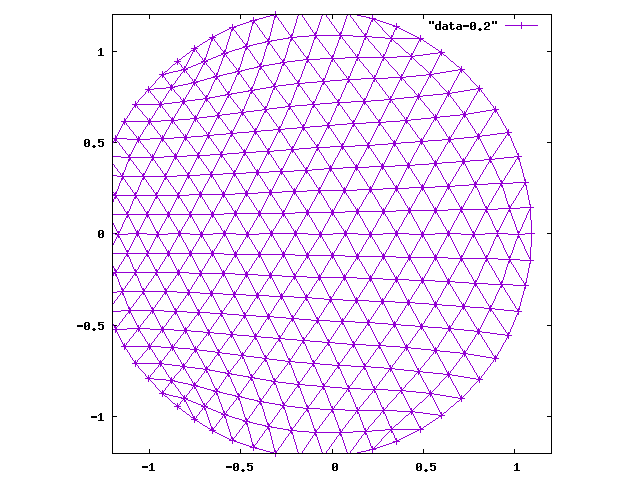}
\includegraphics[width=2.0cm]{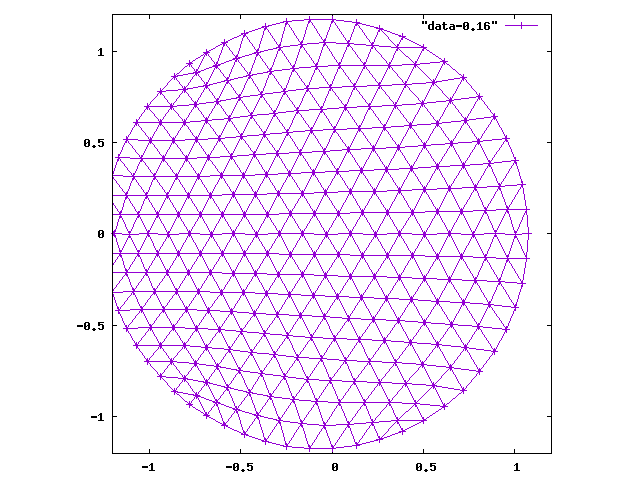}
\includegraphics[width=2.0cm]{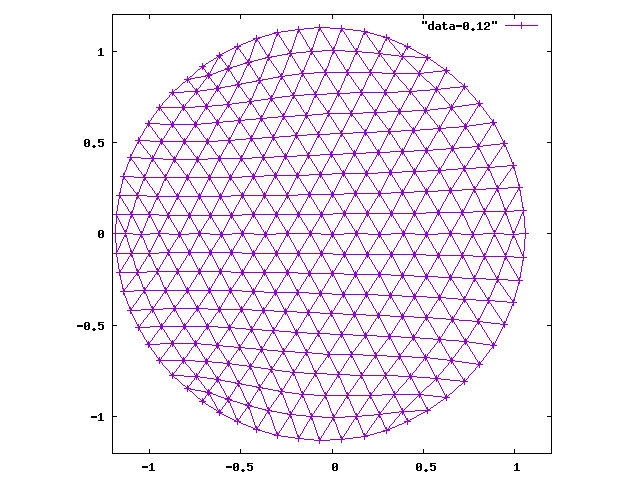}
\includegraphics[width=2.0cm]{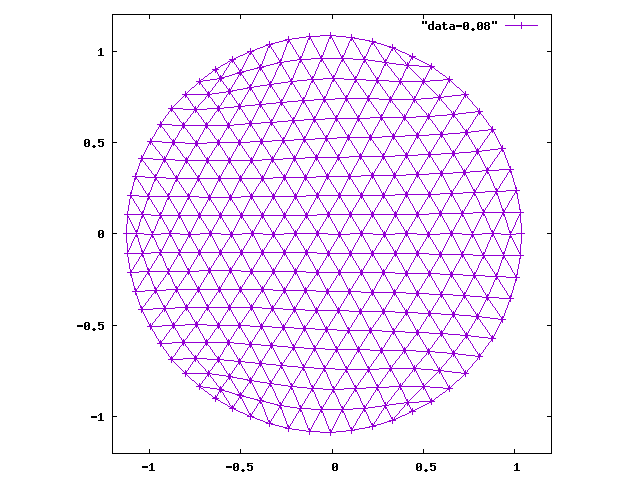}
\includegraphics[width=2.0cm]{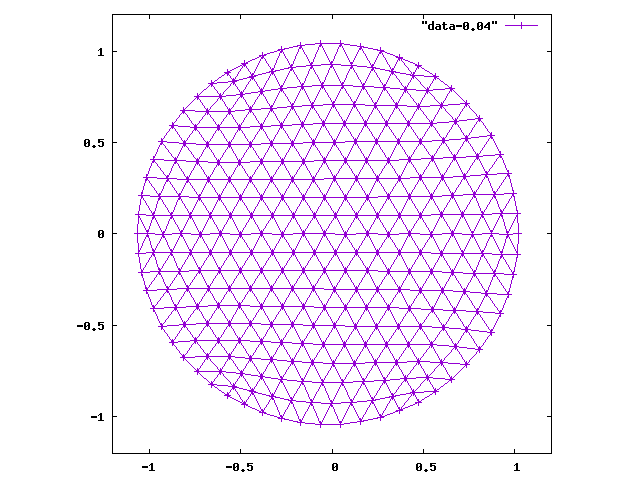}
\includegraphics[width=2.0cm]{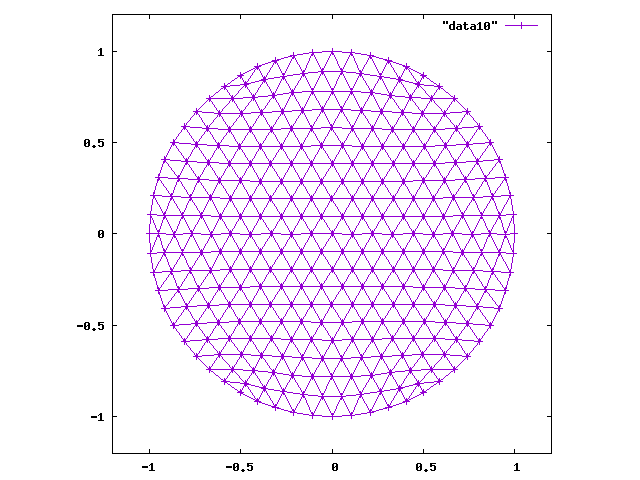} \\
\includegraphics[width=2.0cm]{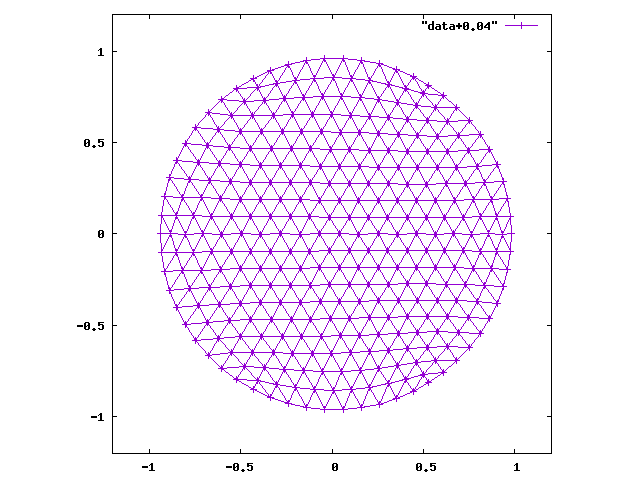}
\includegraphics[width=2.0cm]{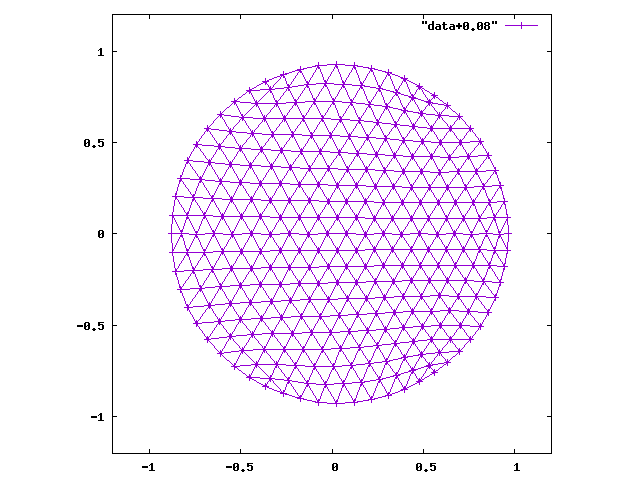}
\includegraphics[width=2.0cm]{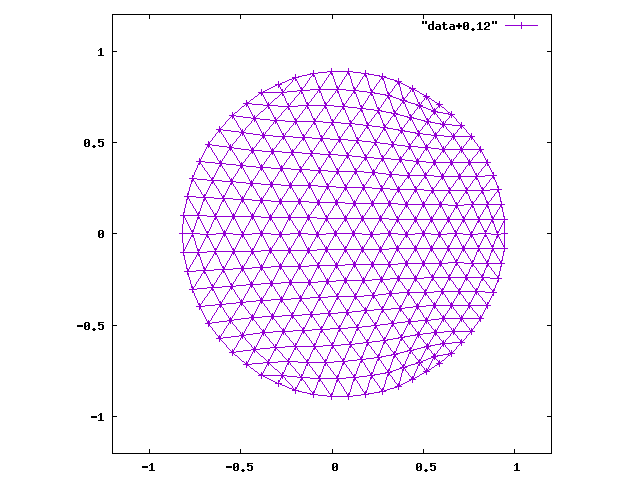}
\includegraphics[width=2.0cm]{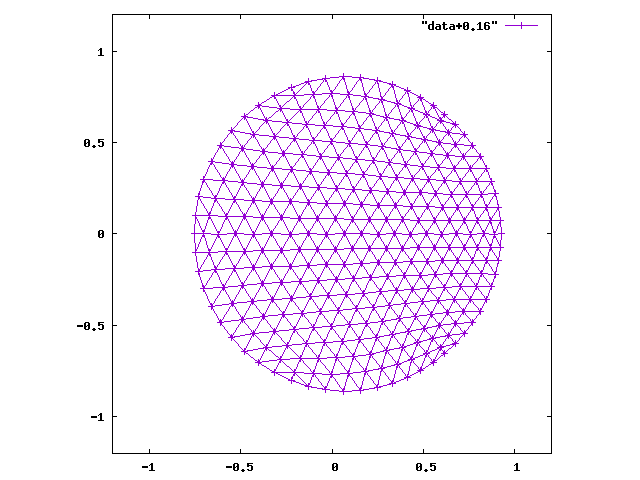}
\includegraphics[width=2.0cm]{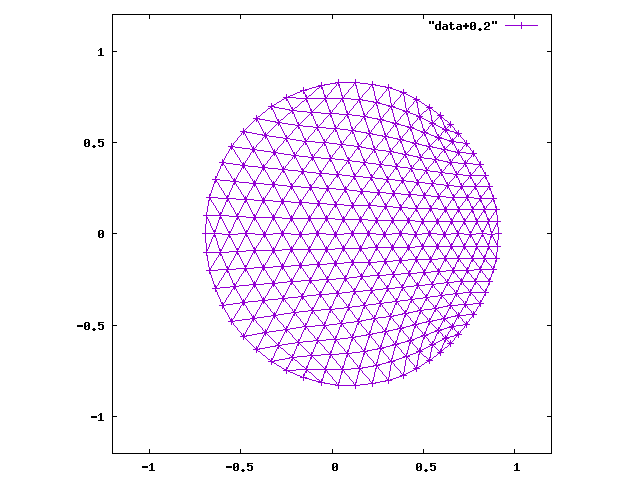}
\caption{The images of the transformation of $D$.
The value of $t$ are from $t=-0.2$  (first row, leftmost) 
to $t = 0.2$ (second row, rightmost) with $0.4$ increments.}
\end{figure}

\begin{figure}[htbp]
\centering
 \includegraphics[width=6.0cm]{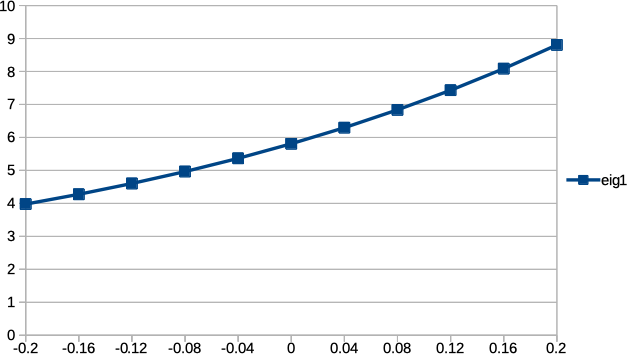} \,
 \includegraphics[width=6.0cm]{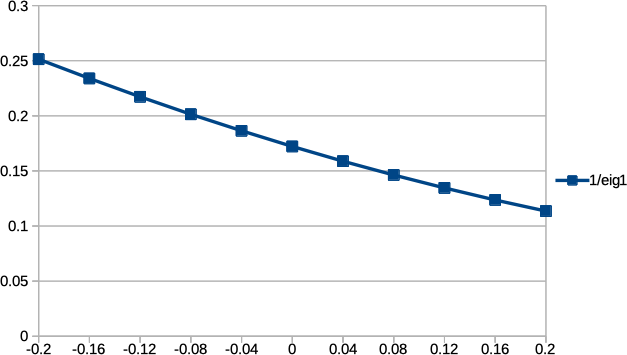}
 \caption{The profiles of the first eigen value (left) and its
reciprocal (right).}
 \label{fig22}
\end{figure}

Then, the eigenvalues of 
\begin{align}
  - \Delta u = \lambda_k(t) u \quad \text{ in } \; g_t(D), \qquad
   u = 0 \quad \text{ on } \; \partial g_t(D)
   \label{eq60}
\end{align}
are computed by the piecewise linear finite element method.
By the calculus in the proof of Theorem~\ref{thm8} and
Theorem~\ref{thm7}, we have
\begin{align*}
  \frac{d^2 a_t}{dt^2} = 2|\sin z - 1|^2 \ge 0 \qquad
 \text{ and } \qquad
  \frac{d^2\;}{dt^2} \frac{1}{\lambda_1(t)} \ge 0,
\end{align*}
although $\cos z$ is not univalent
on $D$.  The computed profiles of the first
eigenvalue and its reciprocal given in Figure~\ref{fig22}
are clearly consistent with Theorem~\ref{thm7}.
\vspace{4mm}

\begin{figure}[htbp]
\centering
 \includegraphics[width=6.0cm]{circle2-5.png} \;
 \includegraphics[width=6.0cm]{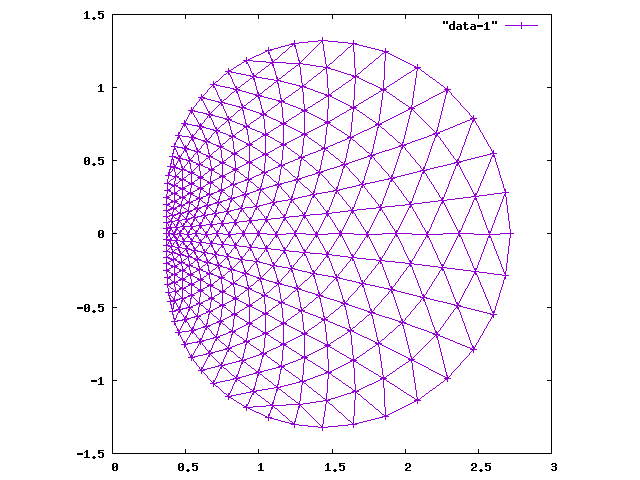} 
 \caption{The unit disk  and its image by $g_1$.}
 \label{fig3}
\end{figure}

Next, we define
\begin{align*}
  g_t(z) = (1 - t)z + t e^z, \qquad z \in \C
\end{align*}
to confirm the statement of Theorem~\ref{thm8}. Note that
$e^z$ is univalent on $D$.  Then, the images of $D = g_0(D)$ and
$g_1(D)$ are given in Figure~\ref{fig3}.
\vspace{4mm}
The first five eigenvalues of \eqref{eq60} on $D$ and $g_1(D)$ are
given in Table~1. Because $a_1 = 1$ in this case, the numbers on
$\lambda_1$ are consistent with the inequality \eqref{37}.
Although Theorem~\ref{thm8} and Remark~4 claim only on the first
eigenvalue $\lambda_1$, we observe similar inequalities hold
for $\lambda_k$ ($2 \le k \le 5$).  

\begin{table}[htbp]
\caption{The first five eigenvalues of Laplacian on $D$ and $g_1(D)$.}
\begin{center}
\begin{tabular}{|c|c|c|c|c|c|}
\hline
   & $\lambda_1$  & $\lambda_2$ & $\lambda_3$ & $\lambda_4$ &
  $\lambda_5$ \\
\hline 
$D$ & 5.80728 & 14.8489 & 14.8489 & 26.9304 & 26.9304 \\
\hline
$g_1(D)$ & 3.69736 & 8.96092 & 10.0331 & 16.8943 & 17.2069 \\
\hline
\end{tabular}
\end{center}
\label{table1}
\end{table}

%At last, we consider the transformation
%\begin{align*}
%  T_t(z) = z + a_2 tz^2 + a_3 t^2z^3 + a_4 t^3 z^4 +a_5 t^4z^5,
%\end{align*}
%where $a_2 = 1 + 0.5i$, $a_3 = -a_2$, $a_4 = 1+i$, $a_5 = - a_4$.
%Although this transformation does not satisfy the condition
%\begin{align}
%   \frac{d^2 a_t}{dt^2} \ge 0,
%   \label{eq61}
%\end{align}
%the reciprocal of the first eigenvalue is still convex as is
%shown in Figure~\ref{fig4}.
%Therefore, we infer that the condition \eqref{eq61} is sufficient but
%not necessary for the convexity of the first eigenvalue of
%Laplacian.
%\vspace{4mm}
%\begin{figure}[htbp]
%\centering
% \includegraphics[width=5.8cm]{ex2-2.png} \,
% \includegraphics[width=5.8cm]{ex2-3.png}
% \caption{The profiles of the first eigen value (left) and its
%reciprocal (right) for $-0.2 \le t \le 0.2$}
% \label{fig4}
%\end{figure}

{}
\end{document}